\documentclass[12pt,leqno]{amsart}
\usepackage[dvips]{graphicx}
\usepackage{amsfonts}
\usepackage{amssymb}
\usepackage{amsmath}
\usepackage{amscd}
\usepackage{graphicx} 
\def\DATE{\today}
\newtheorem{theorem}{Theorem}
\newtheorem{definition}[theorem]{Definition}

\newtheorem{lemma}[theorem]{Lemma}

\newtheorem{proposition}[theorem]{Proposition}

\newcommand\ld{\ldots}

\newcommand\N{\mathbb{N}}
\newcommand\g{\mathfrak{g}}

\newcommand\h{\mathfrak{h}}

\newcommand\K{\mathbb{K}}

\newcommand\p{\mathcal{P}}

\newcommand{\im}{{\rm{Im}}}

\newcommand\pf{\noindent{\it Proof. }}
\newcommand\lr{\left\{ \begin{array}{l}}
\def\ds{\displaystyle}

\title{$k$-step nilpotent  Lie algebras }
\author{Michel Goze, Elisabeth Remm}
\date{Nissan 1, 5775}
\address{Universit\'{e} de Haute Alsace, LMIA, 4 rue des Fr\`{e}res Lumi\`{e}%
re, 68093 Mulhouse}
\email{michel.goze@uha.fr, elisabeth.remm@uha.fr}

\begin{document}

\maketitle

\begin{abstract}

The classification of complex of real finite dimensional Lie algebras which are not semi simple is still in its early stages. For example the nilpotent Lie algebras are classified only up to the dimension 7. Moreover, to recognize a given Lie algebra in a classification list is not so easy. In this work we propose a different approach to this problem. We determine families for some fixed invariants, the classification follows by  a deformation process  or contraction process. We focus on the case of 2 and 3-step nilpotent Lie algebras. We  describe in both cases a deformation cohomology for this type of  algebras and  the  algebras which are rigid regarding this cohomology. Other $p$-step nilpotent Lie algebras are obtained by contraction of the rigid ones.
\end{abstract}

\section{Introduction}\label{sII}

 Let $\g$ be a Lie algebra over an algebraically closed field $\K$ of characteristic $0$. We denote by $\left\{\g^k\;|\; k=0, 1,2,\ld\right\}$ the lower central series of $\g$ defined by $\g^0=\g$ and $\g^{k}=[\g^{k-1},\g]$, for $k=1,2,\ld$  In any Lie algebra the lower central series is a \textit{filtration} in the sense that $[\g^i,\g^j]\subset \g^{i+j+1}$. 
One calls $\g$\textit{ nilpotent} if there is an integer $n$ such that $\g^{n}=\{ 0\}$. If moreover $\g^{n-1}\ne\{ 0\}$ then $n$ is called the\textit{ nilpotent index} or\textit{ nilindex} of $\g$ and denoted by $n(\g)$.
 A nilpotent Lie algebra $\g$ with nilindex $n(\g)$ is called $n(\g)$-step nilpotent. From Engel Theorem, each adjoint operator $ad(X)$ is nilpotent with $ad(X)^{n(\g)}=0$ for any $X \in \g.$ In the present work, 
 we regard a  $p$-step nilpotent Lie algebra as a Lie algebra whose Lie bracket satisfies a $p$-associativity condition. We define for $p$ equal to $2$ or $3$ new cohomological spaces whose main property is that the second group of cohomology characterizes the deformations of $p$-step nilpotent Lie algebras in the class of $p$-step nilpotent Lie algebras.
 Recall that if we 
 denote by $H^*_C(\g, \g)$ the Chevalley-Eilenberg cohomological spaces of $\g,$ the space $H^2_C(\g,\g)$ parametrizes the deformations of $\g$ but in general such a deformation is not nilpotent. Cohomological spaces also different from the 
 Chevalley-Eilenberg's ones are defined in \cite{Vergne} and more adapted to nilpotent Lie algebras. In fact,   if $p$ is fixed, the $2$-cocycles determine $p$-step nilpotent cocycles of $\g / \g^{p+1}$ by 
 $\{ X \in \g , ad(Y)^p(X)=0 \ \forall Y \in \g \}.$ Our point of view is a little bit different. We will describe cohomological spaces which are associated to deformations of $p$-step nilpotent Lie algebras in this family of $p$-step nilpotent Lie algebras.
 
 Let $pNilp_n$ be the variety of $n$-dimensional $p$-step nilpotent Lie algebras over $\K$. Recall that the characteristic sequence $c(\g)$ of a nilpotent Lie algebra is the invariant, up to isomorphism, defined as follows (see \cite{Ancochea-Goze,MGoran}): if $X$ is in $\g$ then the linear operator $ad(X)$ is nilpotent and we denote by  $c(X)$ the decreasing sequence of the dimensions of Jordan blocks of $ad(X).$ The characteristic sequence of $\g$ is  
 $$c(\g)= \max \{ c(X), X \in \g-\g^1 \},$$
 considering the lexicographic order. 
If $\g\in pNilp_n,$ then $c(\g)=(p,\cdots, p,p-1, \cdots, 1)$ where $p$ appears $k_p$ times, $p-1$ appears $k_{p-1}$ times, $\cdots$, and $1$ appears $k_1$ times with $k_p\neq 0 \neq k_1$ and $n=\sum_{i=1}^p i k_i.$ It is clear that if $\g_1 \in pNilp_n$ is a formal deformation of $\g \in pNilp_n$ in the Gerstenhaber sense (\cite{Ge,MGoran,Goze-Remm-Valued})  then $c(\g) \leq c(\g_1)$. If the characteristic sequence of the Lie algebra $\g$ is  $c(\g)=(p,\cdots, p, l ,1)$ with $l\equiv n-1 \ [p]$  then any deformation $\g_1 \in pNilp_n$ of $\g$ has the same characteristic sequence. We denote by $\mathcal{F}^{n,p}_{k_p,k_{p-1}, \cdots, k_1}$ the family of $n$-dimensional $p$-step nilpotent Lie algebras with characteristic sequence $(p,\cdots, p,p-1, \cdots, 1)$ where $i$ appears $k_i$ times, for $1\leq i \leq p.$ The previous remark shows that  $\mathcal{F}^{pk_p+l+1,p}_{k_p,k_l=1,k_1=1}$ and  $\mathcal{F}^{pk_p+2,p}_{k_p,k_1=2}$ (the non mentioned $k_i$ are zero) are open sets in $pNilp_n$ ($n=\sum_{i=1}^p ik_i$). The aim of this paper is to study these open sets for $p=2$ and $3.$

The study of these families is interesting, at least, for two reasons:
\begin{enumerate}
  \item The important works of Vergne on nilpotent Lie algebras lead to conjecture that "there exist no rigid nilpotent Lie algebras (in the variety of Lie algebras )", where rigidity means that any deformation is isomorphic. This conjecture has still not be settled.  In the present work, we define a notion of rigidity  on $pNilp_n$ associated with  deformations internal to this variety. If there exists a counterexample to the Vergne's  conjecture, necessarily this algebra will also be rigid in $pNilp_n$ for a good $p$.
  \item The classification of nilpotent Lie algebras is still an open problem. General results are known for the dimensions smaller than  or equal to $7$ \cite{Gong} and partial results for the dimension $8$ (\cite{Zaili}). In the present work we determine rigid Lie algebras in $\mathcal{F}^{n,p}_{k_p,0,  \cdots, 0 , k_1=1} $ and $\mathcal{F}^{n,p}_{k_p,0,  \cdots, 0 , k_1=2} $ for $p\leq 3$. Then the classification of the elements of this family is in part determined by the contractions of the rigid Lie algebras contained in this variety. 
\end{enumerate} 

\smallskip

We give also some examples of study of families $\mathcal{F}^{n,p}_{k_p,k_{p-1}, \cdots, k_1}$ associated to $p$-step nilpotent Lie algebras with a characteristic sequence which is not maximal, in order to show how it can be useful to classify the $pNilp_n.$

\medskip

\noindent{\bf Notations.} If $\g$ is a Lie algebra over $\K$ with Lie bracket  $\mu$, we denote by $H_C^*(\g,\g)$ or $H_C^*(\mu,\mu)$ its Chevalley-Eilenberg cohomology. The corresponding coboundary operator will be denoted $\delta_{C,\mu}$, or when no confusion is possible, $\delta_C$. 
 \section{Associative Lie multiplication and $2$-step nilpotent Lie algebras}

\subsection{The variety $2Nilp_n$}

Let $\g$ be a Lie algebra over an algebraic field $\K$ of characteristic $0$. If we denote by $[X,Y]$ the Lie bracket of $\g$, it satisfies the following identities
$$
\left\{
\begin{array}{l}
\lbrack X,Y]=-[Y,X], \\
\lbrack \lbrack X,Y],Z]+\lbrack \lbrack Y,Z],X]+\lbrack \lbrack Z,X],Y]=0 \ \ {\text{ \rm (Jacobi Identity),}}
\end{array}
\right.
$$
for any $X,Y,Z \in \g$.
We assume moreover that the Lie bracket is also an associative product, that is, it satisfies
$$[[X,Y],Z]=[X,[Y,Z]],$$
for any $X,Y,Z \in \g$.
The Jacobi identity therefore implies
\begin{equation}
\label{2step}
[[Z,X ],Y]= 0.
\end{equation}

\begin{proposition}
The Lie bracket of the Lie algebra $\g$ is an associative product if and only if $\g$ is a $2$-step nilpotent Lie algebra.
\end{proposition}
In fact, the relation $[[Z,X ],Y]= 0$ means that $[[\g,\g],\g]=\mathcal{C}^2(\g)=0$ where $\mathcal{C}^i(\g)$ denotes the ideals of the descending central series of $\g$. The converse is obvious.

Let $n$ be a fixed integer. Any $n$-dimensional $\K$-Lie algebra $\g$ is identified with the pair $(\K^n, \mu)$ where $\mu$ is the Lie bracket on the underlying vector space of $\g,$ that is $\K^n.$ We fix a basis of $\{ e_i \}$ of $\K^n$. Then $\mu$ is determinated by its structure constants $\displaystyle \mu(e_i,e_j)=\sum_{k=1}^n C_{ij}^k e_k.$ Skew-symmetry and Equation (\ref{2step})  are equivalent to 
\begin{equation}
\label{ }
\left\{
\begin{array}{ l}
\medskip 
C_{ij}^k=-C_{ji}^k ,  \\ 
\displaystyle  \sum_{l=1}^n C_{ij}^l C_{lk}^s=0 , 
\end{array}
\right.
\ \forall k,s \in \{ 1, \cdots , n \}, 1\leq i <j \leq n.
\end{equation}
These polynomial identities define an algebraic variety on $\K^{n^3}$ denoted $2Nilp_n$ in which any point is a $2$-step $n$-dimensional nilpotent Lie algebra over $\K.$

This variety is fibered by the orbits of the action of $GL(n,\K)$ on $2Nilp_n$ 
$$GL(n,\K)\times 2Nilp_n \rightarrow 2Nilp_n$$
given by $(f,\mu)\mapsto f \cdot \mu$ with $f \cdot \mu (X,Y)=f^{-1}\mu(f(X),f(Y)).$ The orbit $\theta(\mu)$ of $\mu$ is the set of Lie algebras isomorphic  to $\g=(\K^n, \mu).$

\begin{definition}
A Lie algebra $\mu \in 2Nilp_n$ is called rigid in $2Nilp_n$ if its orbit $\theta(\mu)$ is open in $2Nilp_n.$
\end{definition}

This notion of rigidity is equivalent to say that any formal deformation $\mu_t=\mu + \sum t^i \phi_i$ of $\mu$ in $2Nilp_n$ is isomorphic to $\mu$ (For general definition of formal deformations see \cite{Ge,Markl,MGoran}). We know that there exists a cohomological sequence $H^*_{def}(\mu, \mu)$ which parametrizes  the deformations of $\mu.$ We shall define this cohomology.

\subsection{The CH-cohomology of $2$-step nilpotent Lie algebras}

Let $\g_0=(\K^n,\mu_0)$ be a  $2$-step nilpotent Lie algebra over $\K$. The CH-cohomology of $\g_0$ is the cohomology associated with the complex 
$(\mathcal{C}^i(\K^n,\K^n), \delta^i_{CH,\mu_0})$ where 
$\mathcal{C}^i(\K^n,\K^n)$ is the vector space of skew-symmetric  $i$-linear maps on $\g_0$ with values in $\g_0$ and 
$\delta^i_{CH,\mu_0}:\mathcal{C}^i(\K^n,\K^n) \rightarrow \mathcal{C}^{i+1}(\K^n,\K^n)$ is the linear operator defined by
$$\left\{
\begin{array}{ll}
\medskip
 \delta^{2p}_{CH,\mu_0}\psi(X_1, \cdots,X_{2p+1})= & \mu_0(X_1, \psi(X_2, \cdots,X_{2p+1}))\\
\medskip
&\displaystyle + \sum_{i=1}^p
\psi(X_1, \cdots,\mu_0(X_{2i},X_{2i+1}), \cdots ,X_{2p+1}),\\
\medskip
\delta^{2p-1}_{CH,\mu_0}\psi(X_1, \cdots,X_{2p})= & \mu_0(X_1, \psi(X_2, \cdots,X_{2p}))\\
 &\displaystyle + \sum_{i=1}^{p-1}
\psi(X_1, \cdots,\mu_0(X_{2i+1},X_{2i+2}), \cdots ,X_{2p}).
\end{array}
\right.$$
Since $\delta^i_{CH,\mu_0} \circ \delta^{i-1}_{CH,\mu_0} =0$, we define
$$H^i_{CH}(\mu_0,\mu_0)= \ds \frac{Z^i_{CH}(\mu_0,\mu_0)}{B^i_{CH}(\mu_0,\mu_0)}$$
with $Z^i_{CH}(\mu_0,\mu_0)=\{\varphi \in \mathcal{C}^i(\K^n,\K^n), \ \delta^{i}_{CH,\mu_0}\varphi=0\}$ and $B^i_{CH}(\mu_0,\mu_0)=\im 
\delta^{i-1}_{CH,\mu_0}.$

\medskip

We shall prove that this cohomology is the  cohomology of deformations for $2$-step nilpotent Lie algebras. For this, we introduce a quadratic operad with the property, that any algebra over this operad is a $2$-step nilpotent Lie algebra.

An operad is a sequence $\p=\{\p(n), n \in \N^*\}$ of $\K[\Sigma _n]$-modules, where $\K[\Sigma _n]$ is the group algebra associated with the symmetric group $\Sigma_n$, with $comp_i$-operations (see  \cite{M.S.S}). The main example corresponds to the free operad $\Gamma(E)=\{\Gamma(E)(n)\}$ generated by a $\K[\Sigma_2]$-module $E$. An operad $\p$ is called binary quadratic if there is a $\K[\Sigma_2]$-module $E$ and a $\K[\Sigma_3]$-submodule $R$ of $\Gamma(E)(3)$ such that $\p$ is isomorphic to $\Gamma(E)/\mathcal{R}$ where $\mathcal{R}$ is the operadic ideal generated by $\mathcal{R}(3)=R.$
\begin{proposition}There exists a binary quadratic operad, denoted by $2\mathcal{N}ilp$, with the property that any $2\mathcal{N}ilp$-algebra is a $2$-step nilpotent Lie algebra. Moreover, this operad is Koszul.
 \end{proposition}
 In fact, we consider $E=sgn_2,$ that is, the representation of $\Sigma_2$ by the signature, then $\Gamma(E)(3)=sgn_3 \oplus V_2$ where $V_2=\{(x,y,z) \in \K^3, x+y+z=0\}$. Let $R$ be the submodule of $\Gamma(E)(3)$ generated by the vectors 
 $(x_i\cdot x_j)\cdot x_k$, with $i,j,k$ all different. We deduce that $2\mathcal{N}ilp(2)$ is the $\K[\Sigma_2]$-module generated by $x_1\cdot x_2$ with the relation $x_2\cdot x_1=-x_1\cdot x_2$ so it is a $1$-dimensional vector space and $2\mathcal{N}ilp(3)=\{0\}$. Let us prove that this operad satisfies the Koszul property (see \cite{M.S.S}). 
The generating function of a binary quadratic operad $2\mathcal{N}ilp$ is
$$g_{2\mathcal{N}ilp}(x)=\sum_{a \geq 1} \frac{1}{a!}\dim (2\mathcal{N}ilp(a))x^a=x+\frac{x^2}{2}.$$
The dual operad  $2\mathcal{N}ilp^!$ of the operad $2\mathcal{N}ilp$ is the quadratic operad $2\mathcal{N}ilp^! := \Gamma (E^\vee)/(R^\perp)$,
where $R^\perp \subset \Gamma (E^\vee)(3)$ is the annihilator of $R \subset \Gamma(E)(3)$. Thus  the dual operad $(2\mathcal{N}ilp)^!$ is $\Gamma(1\! \! 1)$ the free operad generated by a commutative operation. So
$$\dim (2\mathcal{N}ilp)^!(1)=1, \ \dim (2\mathcal{N}ilp)^!(2)=1, \ \dim (2\mathcal{N}ilp)^!(3)=3, \ \dim (2\mathcal{N}ilp)^!(4)=15$$
and more generally, if we denote by $d_k$ the dimension of $(2\mathcal{N}ilp)^!(k)$, we have
$$
\left\{
\begin{array}{l}
\medskip
\displaystyle d_{2k+1}=\sum_{i=1}^k C_{2k+1}^i d_i d_{2k+1-i},\\
\displaystyle d_{2k}=\sum_{i=1}^{k-1} C_{2k}^i d_i d_{2k-i}+\frac{1}{2}C_{2k}^kd_k^2.
\end{array}
\right.
$$
So the generating function of $2\mathcal{N}ilp^!$ is
$$\displaystyle \sum_{k \geq 1}\frac{d_k}{k!}x^k.$$
If an operad $\p$ is Koszul, then its dual $\p^!$ is also Koszul   and the generating functions are related by the functional equation
$$g_{\p}(-g_{\p ^!}(-x))=x.$$ It is known that $\Gamma(1\! \! 1)$ is Koszul, so also $2\mathcal{N}ilp $ and this implies the proposition.
We can verify that the generating function $g_{2\mathcal{N}ilp}$ of the operad $2\mathcal{N}ilp$ satisfies the  functional equation
$$g_{2\mathcal{N}ilp}(-g_{2\mathcal{N}ilp ^!}(-x))=x.$$
The operadic cohomology of a quadratic operad is described in \cite{M.S.S}. We find the CH-cohomology. Since $2\mathcal{N}ilp$ satisfies the Koszul property, this cohomology coincides with the deformation cohomology for $2\mathcal{N}ilp$-algebras, that is the cohomology which parametrizes the formal deformations of a $2$-step nilpotent Lie algebra.

\noindent{\bf Remark.}
 Let us consider an associative algebra $(A,\cdot)$ where $x \cdot y$ denotes the multiplication in $A$. Thus
$$[x,y]=x\cdot y - y\cdot x$$
is a Lie bracket. This Lie bracket is associative if and only if the multiplication of $A$ satisfies
$$(x \cdot y)\cdot z - (y \cdot x)\cdot z -(z \cdot x)\cdot y + (z \cdot y)\cdot x =0.$$
This kind of nonassociative algebras belongs to the family of nonassociative algebras described in \cite{Goze-Remm-Nonass}. In this work, we define ``natural" non-associative algebras, denoted $G$-associative. Thus we can look what happen when we consider a Lie algebra whose Lie bracket satisfies such non-associative identities. But, we can look very quickly that if a Lie bracket is $G$-associative, then it is also associative and we return to the initial case.

\subsection{Deformations and Rigidity  in $2Nilp_{n}$}

Let $\g_0= (\K^n,\mu_0)$ be a $n$-dimensional $2$-step nilpotent Lie algebra over $\K$.  
Let $\g=(\K^n,\mu)$ be a deformation of $\g_0$ in $2Nilp_n$, that is
$$\mu=\mu_0 + \sum t^i \varphi_i$$
where $\varphi_i$ are skew-symmetric bilinear maps on $\K^n$. The Jacobi condition related to $\mu$ gives in particular
that  $\varphi_1 \in Z_{CH}^2 (\g_0,\g_0)$. A deformation of type $\mu=\mu_0 +t \varphi$ will be called a {\it linear deformation} of $\mu_0$. In this case, since we assume that $\g$ is $2$-step nilpotent,  the Jacobi condition related to $\mu$ is equivalent to 
\begin{equation}
\label{2nilp}
\left\{
\begin{array}{l}
   \delta_{CH, \mu_0}\varphi_1=0 ,    \\
     \varphi_1 \circ_1 \varphi_1=0,
\end{array}
\right.
\end{equation}
where $\circ_1$ is the comp$_1$-operation.  

Since $H^*_{CH}(\g_0,\g_0)$ parametrizes the deformations of $\g_0$ in $2Nilp_n$, Nijenhuis-Richardson Theorem (\cite{NR})  gives:

\begin{theorem}
Let $\g_0$ be a $n$-dimensional $2$-step nilpotent Lie algebra on $\K$. If $H^2_{CH}(\g_0,\g_0)=0$,  then $\g_0$ is rigid in $2Nilp_n$.
\end{theorem}

\begin{proposition} \cite{GR-DGA} 
The $(2p+1)$-dimensional Heisenberg algebra $\h_{2p+1}$ is rigid in $2Nilp_{2p+1}$.
\end{proposition}
\pf The dimension of the algebra of derivations of  $\h_{2p+1}$ is computed in \cite{GP}. We deduce $\dim B^2_{CH}(\h_{2p+1},\h_{2p+1})=p(2p+1).$ Let $\{X_1,\cdots,X_{2p+1}\}$ be the basis of $\h_{2p+1}$ satisfying
$$[X_1,X_2]=\cdots=[X_{2i-1},X_{2i}]=\cdots=[X_{2p-1},X_{2p}]=X_{2p+1},$$
the other brackets are equal to zero or deduced by skew-symmetry. 
If $\ds \varphi(X_i,X_j)=\sum_{k=1}^{2p+1}a_{ij}^k X_k$,  then
$\delta_{CH,\mu_0}(\varphi)(X_i,X_j,X_k)=0$
is equivalent to
$$\varphi(\mu_0(X_i,X_j),X_k)+\mu_0(\varphi(X_i,X_j),X_k)=0,$$ that is,
$$
\left\{
 \begin{array}{ll}
\medskip
\displaystyle \varphi(X_1,X_2)=\sum_{k=1}^{2p+1}a_{12}^k X_k;  \\
\medskip
\varphi(X_1,X_i)=a_{1i}^{2p+1} X_{2p+1}  , 3 \leq  i \leq 2p; \ \  \varphi(X_1,X_{2p+1})=-a_{12}^{1} X_{2p+1};\\
\medskip
\varphi(X_2,X_i)=a_{2i}^{2p+1} X_{2p+1}, \  3 \leq  i \leq 2p; \ \   \varphi(X_2,X_{2p+1})=a_{12}^{1} X_{2p+1}; \\
\cdots \\
\medskip
\displaystyle \varphi(X_{2i-1},X_{2i})=\sum_{k=1}^{2p}a_{12}^k X_k+a_{2i-1 \ 2i}^{2p+1}X_{2p+1}, \  2 \leq  i \leq p; \\
\medskip
 \varphi(X_l,X_{s})=a_{ls}^{2p+1} X_{2p+1}, (l,s) \neq (2i-1,2i);\\
\medskip
\varphi(X_{2l},X_{2p+1})=a_{12}^{2l-1} X_{2p+1}, \ l \leq p; \\
\medskip
\varphi(X_{2l-1},X_{2p+1})=-a_{12}^{2l} X_{2p+1} , \ l \leq p. 
\end{array}
\right.
$$
We deduce that $\dim Z^2_{CH}(\h_{2p+1},\h_{2p+1})=p(2p+1).$
So $\dim H^2_{CH}(\h_{2p+1},\h_{2p+1})=0$ and $\h_{2p+1}$ is rigid in $2Nilp_{2p+1}.$

\medskip

Let us note that $\h_{2p+1}$ is not rigid in the variety of $(2p+1)$-dimensional Lie algebras. From \cite{GR-DGA}, any deformation of $\h_{2p+1}$ is a contact Lie algebra. The rigidity in $2Nilp_n$ can be understood by the fact that any nilpotent Lie algebra whose characteristic sequence is $(2,1,\cdots,1)$ is isomorphic to a direct product of an Heisenberg algebra by an abelian algebra. 


\subsection{$2$-step nilpotent Lie algebras with characteristic sequence $(2,\cdots,2,1)$}

Let $\g_{p,1}$ be the $(2p+1)$-dimensional Lie algebras defined by the following brackets given in the basis $\{ X_1 ,\cdots , X_{2p+1} \}$ by
$$\left[X_1,X_{2i}\right]=X_{2i+1}, \ \ 1 \leq i \leq p,$$
and the other brackets are equal to zero or deduced by skew-symmetry. Its characteristic sequence is $(2,\cdots,2,1)$ where $2$ appears $p$ times.

\begin{lemma} \label{lemma1}
Any  $(2p+1)$-dimensional $2$-step nilpotent Lie algebra with characteristic sequence $(2,\cdots,2,1)$ is isomorphic to a linear deformation of $\g_{p,1}.$ 
\end{lemma}
 
\pf Let $\g$ be a $(2p+1)$-dimensional $2$-step nilpotent Lie algebra with characteristic sequence $(2,\cdots,2,1).$ There exists a basis $\{ X_1 ,\cdots ,X_{2p+1} \}$ such that the characteristic sequence is given by the operator $ad(X_1).$  If  $\{ X_1 ,\cdots ,X_{2p+1} \}$ is the Jordan basis of $ad(X_1)$ then the brackets of $\g$ write
$$\left\{
\begin{array}{l}
\medskip
[X_1,X_{2i}]=X_{2i+1}, \ \ 1 \leq i \leq p, \\
\left[X_{2i},X_{2j}\right]=\displaystyle\sum_{k=1}^pa_{2i,2j}^{2k+1}X_{2k+1}, \ \ 1\leq i<j \leq p.
\end{array}
\right.
$$
The change of basis $Y_1=X_1, \ Y_i=tX_i$ for $2\leq i \leq 2p+1$ shows that $\g$ is isomorphic to $\g_t$ whose brackets is
$\mu_t=\mu_0 +t \varphi$ where  $\mu_0$ the multiplication of $\g_{p,1}$ and  $ \varphi(X_{2i},X_{2j})=\sum_{k=1}^pa_{2i,2j}^{2k+1}X_{2k+1}, \ \ 1\leq i<j \leq p$,  $\varphi(X_l,X_s)=0$ for all the other cases with $l<s.$ So $\g$ is a linear deformation of $\g_{p,1}.$

\smallskip

Let us compute the second cohomological group $H_{CH}^2(\g_{p,1},\g_{p,1}).$ It is clear that every $2$-cocycle $\varphi$ is cohomologous to a cocycle $\varphi'$ with $\varphi'(X_1,  X_{2i})=0$ for $1\leq i
\leq p.$  Then we can assume that $\varphi$ satisfies  $\varphi( X_1,Y)=0$ for any $Y$ in $\g_{p,1}.$  In this case, we have 
$$ 0=\delta_{CH,\mu_0}\varphi(X_1,X_{2i},X_l)=\varphi(X_{2i+1},X_l),$$
for any $1 \leq l \leq 2p+1$ and $ 1 \leq i \leq p.$ If we put $\varphi(X_{2i},X_{2j})=\sum_{k=1}^{2p+1} a_{2i,2j}^kX_k,$  the equations $\delta_{CH,\mu_0} \varphi(X_{2i},X_{2j},X_l)=0$ for $l=1,2$ imply that $a_{2i,2j}^{2k}=a_{2i,2j}^{1}=0$ for $1\leq k\leq p.$ But $\delta_{CH,\mu_0} f(X_{2i},X_{2j})=a_{1,2i}X_{2j+1}-a_{1,2j}X_{2i+1}$ and $\varphi$ is cohomologeous to 
\begin{equation}
\label{221}
\left\{
\begin{array}{l}
\medskip
\varphi (X_2,X_4)=\displaystyle \sum_{k=3}^{p}a_{24}^{2k+1}X_{2k+1} ,\\
\medskip
\varphi (X_2,X_{2i})=\displaystyle \sum_{k=2}^{p}a_{2,2i}^{2k+1}X_{2k+1}, \ 3\leq i\leq p, \\
\varphi (X_{2i},X_{2j})=\displaystyle \sum_{k=1}^{p}a_{2i,2j}^{2k+1}X_{2k+1} \ 2\leq i<j \leq p, 
\end{array}
\right.
\end{equation}
if $p \geq 3.$ If $p=2$ the  cocyle $\varphi$  is trivial.  We deduce that for any $2$-step nilpotent $(2p+1)$-dimensional Lie algebra $\g$ with characteristic sequence $(2,\cdots,2,1)$, there is a  $\varphi \in Z_{CH}^2(\g_{p,1},\g_{p,1})$ such that $\g$ is isomorphic to the following Lie algebra:
\begin{equation}
\label{F}
\left\{
\begin{array}{l}
\medskip
\left[X_1,X_{2i}\right]=X_{2i+1}, \ 1 \leq i\leq p, \\
\medskip
\left[X_{2i},X_{2j}\right]=\displaystyle \varphi(X_{2i},X_{2j}), \ 2\leq i<j \leq p.
\end{array}
\right.
\end{equation}
where $\varphi$ satisfies (\ref{221}).  In fact, if $\varphi$ satisfies (\ref{221}), it also satisfies $\varphi \bullet \varphi =0,$ that is,  the Jacobi identity. In particular the subspace $\frak{m}$  generated by $\left\{ X_2, \cdots, X_{2p+1}\right\}$ is a $2$-step or $1$-step nilpotent Lie subalgebra of $\g$. For example, if $p=3$, $\frak{m}$ is a $
6$-dimensional $2$-step nilpotent Lie algebra isomorphic to one of the following:
\begin{itemize}
\item $ \left[X_2,X_{4}\right]=X_7, \left[X_2,X_{6}\right]=X_5, \left[X_4,X_{6}\right]=X_3,$  if $\ c(\frak{m})=(2,2,1,1),$
\item $\left[X_2,X_{6}\right]=X_5,  \left[X_4,X_{6}\right]=X_3$,  if $c(\frak{m})=(2,2,1,1),$
\item $\left[X_2,X_{4}\right]=X_7$ if $c(\frak{m})=(2,1,1,1,1),$
\item $\left[X_2,X_{6  }\right]=X_7$ if $c(\frak{m})=(2,1,1,1,1),$
\item  $\left[X_2,X_{i}\right]=0$ if $\frak{m}$ is abelian.
\end{itemize}
We obtain the classification of $7$-dimensional $2$-step nilpotent Lie algebras with characteristic sequence $(2,2,2,1)$.

 Let $\mathcal{F}^{2p+1,2}_{k_2=p,k_1=1}=\mathcal{F}^{2p+1,2}_{p,1}$ be the family of $2$-step nilpotent $(2p+1)$-dimensional Lie algebras with characteristic sequence $(2,\cdots,2,1)$. It is the orbit, associated with the action of the algebraic group $GL(2p+1, \K)$, of the family  $\mathfrak{f}^{2p+1,2}_{p,1}$ of Lie algebras given in (\ref{F}). Since any $\varphi$ of the family (\ref{F}) is given by (\ref{221}) and satisfies $\varphi \bullet \varphi =0$, the family $\mathfrak{f}^{2p+1,2}_{p,1}$ is parametrized by the structure constants $(a_{ij}^k)$ of $\varphi$ and $\mathfrak{f}^{2p+1,2}_{p,1}$ is a linear plane in $\K^{n^3}$. Then  $\mathcal{F}^{2p+1,2}_{p,1}$  is a connected algebraic variety which is fibered by the orbits of the Lie algebras of $\mathfrak{f}^{2p+1,2}_{p,1}$. In particular, if there exists in $\mathcal{F}^{2p+1,2}_{p,1}$ a rigid Lie algebra in $2Nilp_{2p+1}$, that is, a Lie algebra of $\mathcal{F}^{2p+1,2}_{p,1}$ whose orbit is Zariski-open, its Zariski closure is an algebraic component of $\mathcal{F}^{2p+1,2}_{p,1}$ and coincides with it.  This implies
 \begin{proposition} If $\g$ and $\g'$ are two $(2p+1)$-dimensional Lie algebras rigid in $2Nilp_{2p+1}$ and belonging to $\mathcal{F}^{2p+1,2}_{p,1}$ then they are isomorphic.
 \end{proposition}
 
 Let us determine these Lie algebras. 
 
 \begin{proposition}
\begin{enumerate}
\item The Lie algebra $\g_5=\g_{2,1}$ is rigid in $2Nilp_5.$
  \item Let $\g_7$ be the $7$-dimensional $2$-step nilpotent Lie algebra defined by 
$$[X_1,X_{2i}]=X_{2i+1},  1 \leq i \leq 3, \ \left[X_2,X_{4}\right]=X_7, \ \left[X_2,X_{6}\right]=X_5, \  \left[X_4,X_{6}\right]=X_3.$$
This Lie algebra is rigid in  
 $2Nilp_7$ and any Lie algebra of $\mathcal{F}^{7,2}_{3,1}$ is isomorphic to a contraction of $\g_7$.
  \item Let $\g_9$ be the $9$-dimensional Lie algebra given by
$$[X_1,X_{2i}]=X_{2i+1}, 1 \leq i \leq 4, \  \left[X_2,X_{4}\right]=X_7, \ \left[X_2,X_{8}\right]=X_5, 
\ \left[X_4,X_{6}\right]=X_9, \ \left[X_6,X_{8}\right]=X_3.$$
This Lie algebra is rigid in  
 $2Nilp_9$ and any Lie algebra of $\mathcal{F}^{9,2}_{4,1}$ is isomorphic to a contraction of $\g_9$.
\end{enumerate}
\end{proposition}

\noindent{\it Proof.}  In fact, if $\g \in 2Nilp_{2p+1}$ with Lie bracket  $\mu,$ the linear space of $2$-cocycles $\varphi \in Z^2_{CH}(\g,\g)$ satisfying Equations (\ref{221}) is of dimension $m_p=\ds \frac{p(p+1)(p-2)}{2}$. For $p=2$, we obtain $m_2=0$ and $\g_{2,1}$ is rigid in $2Nilp_5$. For $p=3$ and $4$, we have respectively $m_3= 6$ and $m_4=20$. This corresponds to the dimension of the $2$-coboundaries satisfying $\delta_{CH,\mu}f(X_1,X_i)=0$. We deduce that, in each case, $\dim H^2_{CH}(\g,\g)=0$ and $\g_7$ and $\g_9$ are rigid.

\noindent{\bf Remark.}
Since $\mathcal{F}^{7,2}_{3,1}$ and $\mathcal{F}^{9,2}_{4,1}$ are  connected, they are respectively the closure of the orbit of $\g_7$ and $\g_9$. Then any rigid Lie algebra in $\mathcal{F}^{7,2}_{3,1}$ (respectively $\mathcal{F}^{9,2}_{4,1}$) is isomorphic to $\g_7$ (respectively  $\g_9$). Others Lie algebras of these families are in the boundary of the closure of the orbit of $\g_7$ or $\g_9$ and are, by definition, contractions of $\g_7$ or $\g_9$.  (see \cite{MGoran} for definition of contractions).

 \begin{proposition}
 If $p \geq 5$, any $2$-step nilpotent of dimension $2p+1$ 
with characteristic sequence $(2,\cdots,2,1)$ is not rigid in 
 $2Nilp_{2p+1}.$
\end{proposition}
\noindent{\it Proof.} We first remark that a Lie algebra of $\mathcal{F}^{2p+1,2}_{p,1}$ is rigid in $2Nilp_{2p+1}$ if and only if its second cohomological  space $H^2_{CH}$ is trivial. In fact, Nijenhuis-Richardson Theorem implies that a Lie algebra in $2Nilp_{2p+1}$ with a trivial second cohomology space is rigid in  $2Nilp_{2p+1}.$ For the converse we know that if a rigid Lie algebra has a nontrivial cohomology, its  affine scheme associated with the algebraic variety $2Nilp_{2p+1}$ is not reduced. But the family $\mathcal{F}^{2p+1,2}_{p,1}$ is isomorphic to a plane and its associated affine scheme  is reduced. Then here the converse of Nijenhuis-Richardson Theorem is true and we have to prove that $\dim H^2_{CH}(\g,\g) \neq 0.$ To compute $Z^2_{CH}(\g,\g)$, we can consider only the $2$-cocycles satisfying $\varphi(X_1,Y)=0$ for any $Y \in \g$. This implies that $\varphi$ is completely determinate by its constants
$$\varphi(X_{2i},X_{2j})=\sum a_{ij}^kX_{2k+1}.$$
Then this space is of dimension $\frac{p^2(p-1)}{2}.$ Consider now a linear  endomorphism of $\g$ such that $\delta_{CH,\mu}f(X_1,X_i)=0.$ This conditions implies that $f(X_{2i+1})=[f(X_1),X_{2i}]+[X_1,f(X_{2i})]$. If we write
$$\delta_{CH,\mu}f(X_{2i},X_{2j})=\sum b_{ij}^kX_{2k+1}$$
the coefficients $ b_{ij}^k$ are linear combinations of the coefficients $\alpha_{1,1},\alpha_{1,2i},i=1,\cdots,p$ and
$\alpha_{2i,1},\alpha_{2i,2j}$ for $i,j=1,\cdots,p$ where the coefficients $\alpha_{i,j}$ are the coefficients of the matrix of $f$ in the basis $\{X_1, \cdots , X_{2p+1}\}$. Then the $b_{ij}^k$ are linear combinations of $(p+1)^2$ parameters. But we have $\frac{p^2(p-1)}{2}$ coefficients $a_{ij}^k$. Since 
$$\frac{p^2(p-1)}{2} \leq (p+1)^2 \Rightarrow p < 5$$
then $\dim B^2_{CH}(\g,\g) < \dim Z^2_{CH}(\g,\g)$ as soon as $p \geq 5$. 
\medskip

\noindent As soon as $p\geq5$ there are no rigid Lie algebras in $2Nilp_n$. A family with $k$ parameters of Lie algebras of $2Nilp_n$ is called irreducible if two algebras corresponding to different values of the parameters give two non isomorphic  algebras. A irreducible family with $k$ parameters is rigid if any deformation of an algebra of this family is isomorphic to an algebra of this family. This is equivalent to say that the orbit of this family is open. For $p\geq 5$ there exists an irreducible family  with $k$ parameters rigid in $2Nilp_n$ and any algebra in $2Nilp_n$ is a contraction of an element of this family.

\subsection{$2$-step nilpotent Lie algebras with characteristic sequence $(2,\cdots,2,1,1)$}

Let us denote by $\g_{p-1,2}$ the $(2p)$-dimensional Lie algebra given by the brackets
$$[X_1,X_{2i}]=X_{2i+1}, \ 1\leq i \leq p-1,$$
and other brackets are equal to zero or deduced by skew-symmetry.
\begin{lemma}
Any $2$-step nilpotent $(2p)$-dimensional Lie algebra with characteristic sequence $(2,\cdots,2,1,1)$  is isomorphic to a linear deformation of $\g_{p-1,2}$. Its Lie bracket is isomorphic to $\mu=\mu_0 + t \varphi$,  where $\mu_0$ is the Lie bracket of $\g_{p-1,2}$, and $\varphi$ is a skew-bilinear form such that 
\begin{equation}
\label{def}
\left\{
\begin{array}{l} 
\medskip
    \varphi \in Z^2_{CH}(\g_{p-1,2},\g_{p-1,2}), \\
     \varphi \circ_1 \varphi=0.
\end{array}
\right.
\end{equation}
\end{lemma}

\pf The proof is similar to the proof of Lemma  \ref{lemma1}.  Note that, for a general linear deformation $\mu=\mu_0+t\varphi$, the map $\varphi$ is a $2$-cocycle for the Chevalley-Eilenberg cohomology of $\mu_0$ satisfying also the Jacobi condition  (\cite{NR}).  But, since $\mu$ is $2$-step nilpotent, this condition reduces to $\varphi(\varphi(X,Y),Z)=0$ for any $ X,Y,Z.$ 
\begin{lemma}
Any cocycle  in $Z^2_{CH}(\g_{p-1,2},\g_{p-1,2})$ is cohomologous to a cocycle satisfying
\begin{equation}
\label{Z2kk}
\left\{
\begin{array}{l}
\medskip
\varphi(X_1,X_{i})=0,  i=1,\cdots 2p-1, \ \varphi(X_1,X_{2p})=a X_{2p},\\
\medskip
\varphi(X_{2j+1},X_{k})=0,   i=1, \cdots, p, \ j=1,\cdots,p-1, \ k=1,\cdots ,2p, \\
\medskip
\varphi(X_{2i},X_{2j})=\displaystyle \sum_{k=1}^{p-1}a_{2i,2j}^{2k+1}X_{2k+1}+a_{2i,2j}^{2p}X_{2p}, \ 2\leq i<j \leq p.   
\end{array}
\right.
\end{equation}
\end{lemma}
\noindent The proof follows from a direct computation. The coboundaries of Type (\ref{Z2kk}) satisfy
$$\left\{
\begin{array}{l}
\medskip 
\delta f(X_{2i},X_{2j})=a_{1,2i}X_{2j+1}-a_{1,2j}X_{2i+1}, \ 1 \leq i<j \leq p-1,\\
\medskip
\delta f(X_{2i},X_{2p})=-a_{1,2p}X_{2i+1}, \ 1 \leq i<\leq p-1.
\end{array}
\right.
$$
We deduce
\begin{equation}
\label{H2kk }
\dim H^2_{CH}(\g_{p-1,2},\g_{p-1,2})=\ds \frac{p^3-p^2-2p+2}{2}.
\end{equation}
For example, $\dim H^2_{CH}(\g_{p-1,2},\g_{p-1,2})=0$ if and only if $p=1$ and in this case $\g_{0,2}$ is the $2$-dimensional abelian Lie algebra and we know that this Lie algebra is rigid in the variety of nilpotent Lie algebra of dimension $2$. If $p >1$, then $H^2_{CH}(\g_{p-1,2},\g_{p-1,2})$ is not trivial.

Let $\mathcal{F}^{2p,2}_{k_2=p-1,k_1=2}=\mathcal{F}^{2p,2}_{p-1,2}$ be the family of $2$-step nilpotent $(2p)$-dimensional Lie algebra with characteristic sequence $(2,\cdots,2,1,1)$. Any element of this family is isomorphic to a linear deformation of $\g_{p-1,2}$ whose bracket satisfies $\mu=\mu_0 + \varphi$ with $\varphi$ satisfies (\ref{def}). Let $\varphi \in  Z^2_{CH}(\g_{p-1,2},\g_{p-1,2})$ be given by (\ref{Z2kk}). The condition $\varphi(\varphi(X,Y),Z)=0$   is equivalent to:
\begin{enumerate}
  \item $a=0, a_{2i,2p}^{2p}=0, \ 1 \leq i \leq p-1$
  \item and one of the following conditions:
  \begin{enumerate}
  \item $a_{2i,2j}^{2p}=0, 1 \leq i < j \leq p-1$
  \item $\varphi(X_{2i},X_{2p})=0, \  \leq i \leq p-1$
\end{enumerate} 
\end{enumerate}

\begin{proposition}
The family $\mathcal{F}^{2p,2}_{p-1,2}$  of $2$-step nilpotent $(2p)$-dimensional Lie algebras with characteristic sequence $(2,\cdots,2,1,1)$
 is the union of two algebraic components, the first one, $\mathcal{C}_1(\mathcal{F}^{2p,2}_{p-1,2})$, corresponds to the cocycles
\begin{equation}
\label{C1}
\varphi(X_{2i},X_{2j})=\displaystyle \sum_{k=1}^{p-1}a_{2i,2j}^{2k+1}X_{2k+1}, \ 2\leq i<j \leq p,   
\end{equation}
the second one, $\mathcal{C}_2(\mathcal{F}^{2p,2}_{p-1,2})$,  to the cocyles
\begin{equation}
\label{C_2}
\varphi(X_{2i},X_{2j})=\displaystyle \sum_{k=1}^{p-1}a_{2i,2j}^{2k+1}X_{2k+1}+a_{2i,2j}^{2p}X_{2p}, \ 2\leq i<j \leq p-1,   
\end{equation}
where the non defined product $\varphi(X,Y)$ are nul or obtained by skew-symmetry. 
\end{proposition}
Each one of these components is a  regular algebraic variety. These  components can be characterized by a property of the center: it is generated by $\{X_{2i+1}, \ i=1,\cdots, p-1\}$ for $\mathcal{C}_1(\mathcal{F}^{2p,2}_{p-1,2})$ and by  $\{X_{2i+1}, \ i=1,\cdots, p-1 \ , X_{2p}\}$ for $\mathcal{C}_2(\mathcal{F}^{2p,2}_{p-1,2})$.

\begin{proposition} \label{prop16}
\begin{enumerate}
  \item Assume that $p \geq 5$. A $2p$-dimensional Lie algebra  of  $\mathcal{C}_1(\mathcal{F}^{2p,2}_{p-1,2})$ is never rigid. 
  \item Let $p=3$ and $\g_6$ be the Lie algebra of $\mathcal{C}_1(\mathcal{F}^{6,2}_{2,2})$ given by
  $$[X_1,X_{2i}]=X_{2i+1}, i=1,2, \ [X_2,X_6]=X_5, \ [X_4,X_6]=X_3.$$ Then $\g_6$ is rigid in $2Nilp_{6}.$
  \item  Let $p=4$  and $\g_8$ be the Lie algebra of $\mathcal{C}_1(\mathcal{F}^{8,2}_{3,2})$ given by
  $$[X_1,X_{2i}]=X_{2i+1}, i=1,2,3, \ [X_{2},X_{4}]=X_{7},  \ [X_{4},X_{8}]=X_3,  \ [X_{6},X_{8}]=X_5.$$ Then $\g_8$ is rigid in $2Nilp_{8}.$
\end{enumerate}
\end{proposition}

\pf \begin{enumerate}
  \item Let $\g$ be in $\mathcal{C}_1(\mathcal{F}^{2p,2}_{p-1,2})$. A $2$-cocycle belonging to $Z^2_{CH}(\g,\g)$ is cohomologous to a cocycle of type (\ref{C1}). The space of such cocycles is of dimension $\frac{p(p-1)^2}{2}$. The coboundaries belonging to this space are parametrized by the coefficients 
  $a_{11}, a_{1,i}$ and $a_{i,1}$ for  $i=2, \cdots, 2p$, $a_{2k,2l}$ for $k,l=1, \cdots, p,$ that is, $(p+1)^2$ coefficients.
  As soon as $p\geq 5$ we have $  \frac{p(p-1)^2}{2}> (p+1)^2.$ Moreover any rigid Lie algebra in $\mathcal{C}_1(\mathcal{F}^{2p,2}_{p-1,2})$ has a second group of cohomology of dimension $0.$ So there are no rigid Lie algebras  in $\mathcal{C}_1(\mathcal{F}^{2p,2}_{p-1,2})$ for $p \geq 5.$ 
  \item  and (3) The Lie algebras $\g_6$ and $\g_8$ satisfy $\dim H^2_{CH}(\g,\g)=0.$ Thus they are rigid.  Remark that the algebra $\g_8$ corresponds to the algebra $N_9^{8,3}$ in the classification of \cite{Zaili}.  The algebra $N_5^{8,3},N_6^{8,3},N_7^{8,3},N_8^{8,3}$ which are the other algebras belonging to $\mathcal{C}_1(\mathcal{F}^{8,2}_{3,2})$ are contractions of the algebra $N_9^{8,3}.$
\end{enumerate}

\begin{proposition}\begin{enumerate}
  \item Let $p=3$ and $\frak{h}_6$ be the Lie algebra of $\mathcal{C}_2(\mathcal{F}^{6,2}_{2,2})$ given by
  $$[X_1,X_{2i}]=X_{2i+1}, i=1,2, \ [X_2,X_4]=X_6.$$ Then $\frak{h}_6$ is rigid in $2Nilp_{6}.$
     \item Let $p=4$  and $\frak{h}_8$ be the Lie algebra of $\mathcal{C}_2(\mathcal{F}^{8,2}_{3,2})$ given by
  $$[X_1,X_{2i}]=X_{2i+1}, i=1,2,3, \ [X_{2},X_{6}]=X_{5},  \ [X_{2},X_{4}]=X_8.$$ Then $\frak{h}_8$ is rigid in $2Nilp_{8}.$
   \item Let $p=5$  and $\frak{h}_{10}$ be the Lie algebra of $\mathcal{C}_2(\mathcal{F}^{10,2}_{4,2})$ given by
  $$[X_1,X_{2i}]=X_{2i+1}, i=1,2,3,4,  \ [X_{2},X_{4}]=X_{10}, \ [X_{2},X_{6}]=X_{5},  \ [X_{2},X_{8}]=X_3, $$$$[X_{4},X_{6}]=X_{9}, \ [X_{4},X_{8}]=X_{7},  \ [X_{6},X_{8}]=X_3.$$ Then $\frak{h}_{10}$ is rigid in $2Nilp_{10}.$
  \item Assume that $p \geq 6$. A $2p$-dimensional Lie algebra  of  $\mathcal{C}_2(\mathcal{F}^{2p,2}_{p-1,2})$ is never rigid in $2Nilp_{2p}$ and also in the variety $Lie(2n)$ of $2n$-dimensional Lie algebras.
\end{enumerate}
\end{proposition}

The proof of the proposition is similar to the proof of Proposition \ref{prop16}.

\medskip

\section{Cubic associative Lie multiplication}

\subsection{The variety $3Nilp_n$ of $n$-dimensional $3$-step nilpotent Lie algebras and cubic associative algebras}

Let $A$ be a $\K$-associative algebra with binary multiplication $xy$. The associativity which is the quadratic relation
$$(xy)z=x(yz)$$
implies six cubic relations
\begin{equation}
\label{cub}
\lr
((xy)z)t=(x(yz))t,\\
(x(yz))t=x((yz)t),\\
x((yz)t)=x(y(zt)), \\
x(y(zt))=(xy)(zt),\\
(xy)(zt)=((xy)z)t.
\end{array}
\right.
\end{equation}
Recall that these relations correspond to the edges of the Stasheff pentagon.
\begin{definition}
\label{cubic associative}
A binary algebra, that is, an algebra whose multiplication is given by a bilinear map, is called cubic associative if the multiplication satisfies the cubic relations (\ref{cub}).
\end{definition}
We call these relations cubic because if we denote by $\mu$ the multiplication, it occurs exactly three times in each  term of the relations. For example, the first relation can be written
$$\mu\circ(\mu \circ(\mu \otimes Id)\otimes Id)=\mu\circ(\mu \circ(Id \otimes \mu)\otimes Id)$$
which is cubic in $\mu$. It is the same thing for all other relations.

\begin{proposition}
Let $\g$ be a Lie algebra. The Lie bracket is cubic associative if and only if $\g$ is $3$-step nilpotent.
\end{proposition}
In fact, the first identity of (\ref{cub}) becomes
$$
\begin{array}{lll}
[[[X_1,X_2],X_3],X_4] & = & [[X_1,[X_2,X_3]],X_4] \\
&=& -[[[X_2,X_3],X_1],X_4]
\end{array}
$$
and finally
$$[[[X_1,X_2],X_3],X_4]+[[[X_2,X_3],X_1],X_4]=-[[[X_3,X_1],X_2],X_4]=0,$$
which implies that $\g$ is $3$-nilpotent. Conversely, if $\g$ is $3$-nilpotent, all the relations of (\ref{cub}) are satisfied.

\medskip

Let $\g=(\K^n,\mu)$ be a $3$-step nilpotent $n$-dimensional Lie algebra and $\{X_1,\cdots,X_n\}$ be a basis of $\K^n.$ The structure constants $(C_{ij}^k)$ of $\mu$ related to this basis satisfy the following polynomial equations:
\begin{equation}
\label{3s}
\left\{
\begin{array}{l}
\ds C_{ij}^k=-C_{ji}^k,\\
 \ds     \sum_{l=1}^nC_{ij}^lC_{lk}^s +C_{jk}^lC_{li}^s +C_{ki}^lC_{lj}^s =0, \ s=1,\cdots,n,  \\
 \medskip
\ds \sum_{t,u =1}^nC_{ij}^t C_{tk}^u C_{ul}^s =0,   \  s=1,\cdots,n
\end{array}
\right.
\end{equation}
for any $i,j,k,l \in \{1,\cdots,n\}$. 
\medskip
Let $\g$ be a $n$-dimensional $3$-step nilpotent Lie algebra. If we consider the vector space $\K^{n^3}$ parametrized by the structure constants $C_{ij}^k$, the polynomial equations (\ref{3s})  are the equations of the algebraic variety  $3Nilp_n$. Any element of $3Nilp_n$ is  a Lie bracket $\mu$ of a $3$-step nilpotent, $n$-dimensional Lie algebra $\g=(\K^n,\mu)$. We identify $\g$ with its Lie bracket $\mu$.

Let $\g \in 3Nilp_n$. Its characteristic sequence  is  of type $(3,\cdots,3,2,\cdots,2,1,\cdots,1)$ with $n=3k_3+2k_2+k_1$ where $k_i$ is the number of $i$ in the characteristic sequence.  The aim of the next section is to describe the families $\mathcal{F}^{n,3}_{k_3,k_2,k_1}$ of these Lie algebras. For example, in dimension $7$, these families are $\mathcal{F}^{7,3}_{2,0,1}$, $\mathcal{F}^{7,3}_{1,1,2}$ and $\mathcal{F}^{7,3}_{1,0,4}$.

\subsection{Deformations and rigidity in $3Nilp_n$}

Let $\g_{k_3,k_2,k_1}$ be the $n=3k_3+2k_2+k_1$-dimensional nilpotent Lie algebra with characteristic sequence $(3,\cdots,3,2,\cdots,2,1,\cdots,1)$ belonging to $\mathcal{F}^{n,3}_{k_3,k_2,k_1}$ given by
\begin{equation}
\label{mod}
\left\{
\begin{array}{l}
\medskip
[X_1,X_{2+3i}\rbrack=X_{3+3i}, \ [X_1,X_{3+3i}]=X_{4+3i}, \ i=0,\cdots,k_3-1, \\
\medskip
\lbrack X_1,X_{3k_3+2j}\rbrack=X_{3k_3+2j+1}, \ j=1,\cdots k_2, \\
\lbrack X_1,X_{3k_1+2k_2+l+1}\rbrack=0, \ l=1,\cdots,k_1-1.
\end{array}
\right.
\end{equation}

\begin{lemma}
Any  $n=3k_3+2k_2+k_1$-dimensional $3$-step nilpotent Lie algebra $\g,$ i.e $\g \in  \mathcal{F}^{n,3}_{k_3,k_2,k_1},$  is a linear deformation of $\g_{k_3,k_2,k_1}$.
\end{lemma}
The proof is similar to the proof of Lemma \ref{lemma1}.

As consequence, if $\mu_0$ is the Lie bracket of $\g_{k_3,k_2,k_1}$, the Lie bracket $\mu$ of $\g$ is isomorphic to $\mu_0+t\varphi$ where $\varphi$ is a skew-bilinear map satisfying:
\begin{equation}
\label{3step}
\left\{
\begin{array}{l}
\delta_{C,\mu_0}^2(\varphi)=0,\\
\varphi \bullet \varphi =0,\\
\delta_{R,\mu_0}^2(\varphi)=0,\\
(\mu_0 \circ_1\varphi) \circ_1 \varphi + (\varphi \circ_1 \varphi) \circ_1 \mu_0  +(\varphi\circ_1 \mu_0)   \circ_1 \varphi =0, \\
(\varphi  \circ_1\varphi) \circ_1 \varphi =0,
\end{array}
\right.
\end{equation}
where $\delta_{R,\mu_0}^2(\varphi)$ is the $4$-linear map
$$\delta_{R,\mu_0}^2(\varphi)=(\mu_0 \circ_1 \mu_0) \circ_1 \varphi + (\mu_0 \circ_1 \varphi) \circ_1 \mu_0  +(\varphi\circ_1 \mu_0)   \circ_1 \mu_0,$$
the identity
$$\varphi \bullet \varphi =0$$
is the Jacobi identity, and $\circ_1$ is the comp$_1$-operation (composition on the first argument).
If we put $\delta_{R,\mu}^1(f)=\delta_{C,\mu}^1( f)$ for any linear endomorphism $f$ on a Lie algebra then
$$\delta_{R,\mu_0}^2(\delta^1_{R,\mu_0}( f))=0.$$
Let $\g=(\K^n,\mu)$ be a $n$-dimensional $\K$-Lie algebra. 
We consider the complex 
$$(\mathcal{D}^i(\K^n,\K^n), \delta^i_{CR, \mu})_{i\geq 1}$$
 where the space $\mathcal{D}^i(\K^n,\K^n)$ of $i$-cochains is $\mathcal{D}^i(\K^n,\K^n)= \mathcal{C}^i(\K^n, \K^n)  $ for $i=1,2$ and $\mathcal{C}^i(\K^n, \K^n) \oplus \mathcal{C}^{i+1}(\K^n, \K^n)$ for $i\geq 3.$  The coboundary operators are given by $\delta^1_{CR, \mu}=\delta^1_{C, \mu}$, 
$\delta^2_{CR,\mu}(\varphi)=(\delta_{C,\mu}^2(\varphi),\delta_{R,\mu}^2(\varphi)), \  \delta^3_{CR}(\phi,\psi)=(\delta_{C,\mu}^3(\phi),\delta_{R,\mu}^3(\psi))$
and
$$\delta_{R,\mu}^3(\psi)(X_1,\cdots,X_5)=\mu(\psi(X_1,\cdots,X_4),X_5)-\psi(\mu(X_1,X_2),X_3,X_4,X_5).$$
The associated cohomology is the deformation cohomology  for the $3$-step nilpotent Lie algebras. To illustrate, we compute some spaces of this cohomology for particular $3$-step nilpotent Lie algebras and describe the link with the problem of classification.

\medskip

\noindent{\bf Examples}\begin{enumerate}
  \item {\it Characteristic sequence $(3,1,1)$.} Let $\g_{1,0,2}$ be the $5$-dimensional Lie algebra defined by
$$[X_1,X_2]=X_3, \ [X_1,X_3]=X_4.$$
We denote by $\mu_0$ its Lie bracket. Let $\mu_0+t\varphi$ be a $3$-step nilpotent linear deformation of $\mu_0$. The bilinear map $\varphi$ satisfies relation (\ref{3step}). It is cohomologous to a cocycle, still denoted by $\varphi$, which satisfies $\varphi(X_1,Y)=0$ for any $Y.$ This implies that
$Z^2_{CR}(\mu_0,\mu_0)$ is constituted of bilinear forms $\varphi$ defined by
$$\varphi(X_2,X_3)=aX_5, \ \varphi(X_2,X_5)=bX_4+cX_5.$$
The quadratic condition $(\mu_0 \circ_1 \varphi) \circ_1 \varphi + (\varphi \circ_1 \varphi) \circ_1 \mu_0  +(\varphi\circ_1 \mu_0 )  \circ_1 \varphi =0$ implies that 
$ab=0$ and the ternary condition $(\varphi  \circ_1 \varphi) \circ_1 \varphi =0$ gives $c=0$. We deduce that $\varphi$ satisfies
$$\varphi(X_2,X_3)=aX_5, \ \varphi(X_2,X_5)=bX_4,$$
with $ab=0$. We find  the classification of $5$-dimensional  $3$-step nilpotent Lie algebras again.
  \item {\it Characteristic sequence $(3,1,\cdots,1)$ and $n=3+p$}. If $\g$ is such a Lie algebra, its Lie bracket $\mu$ is isomorphic to a linear deformation of $\mu_0$, the Lie bracket of $\g_{1,0, p}$. The bilinear form $\varphi$ satisfies conditions (\ref{3step}), and we can assume that $\varphi(X_1,Y)=0$ for any $Y$. This implies in particular that $\varphi(X_4,Y)=0$. We consider the $2$-cocycles of $Z^2_{CR}(\mu_0,\mu_0)$ satisfying these conditions. They are cohomologous to $2$-cocycles given by
$$
\left\{
\begin{array}{l}
 \varphi(X_2,X_3)=\sum_{i\geq 5}a_{2,3}^iX_i ,\\
\varphi(X_2,X_k)=\sum_{i\geq 4}a_{2,k}^iX_i  , \ 5 \leq k \leq n, \\
 \varphi(X_l,X_k)=\sum_{i\geq 4}a_{l,k}^iX_i  , \  5 \leq l < k  \leq n. \\
\end{array}
\right.
$$
In particular $\dim H^2_{CR}(\mu_0,\mu_0)=\displaystyle (n-3)\frac{n^2-7n+14}{2}-1.$ Since the characteristic sequence is equal to $(3,1,\cdots,1)$ we have necessarily $ \varphi(X_l,X_k)=a_{l,k}^4X_4 $ for $  5 \leq l < k  \leq n$ and also for $l=2, k \geq 5$. If $ \varphi(X_2,X_3)$ is not $0$, then all the constants $a_{l,k}^4$ are $0$ and $\g$ is a direct product of the $5$-dimensional Lie algebra described above and an abelian ideal. If we assume that $\g$ is indecomposable, then $ \varphi(X_2,X_3)=0$ and $\g$ is isomorphic to the Lie algebra
\begin{equation}
\label{clas3111}
\left\{
\begin{array}{l}
[X_1,X_{2}\rbrack=X_{3}, \
\lbrack X_1,X_{3}\rbrack=X_{4},  \\
\lbrack X_2,X_{k}\rbrack=a_{2,k}X_4, \ k \geq 5,\\
\lbrack X_l,X_{k}\rbrack=a_{l,k}^4X_{4}, \ 5 \leq l < k\leq n.

\end{array}
\right.
\end{equation}
\end{enumerate}

\subsection{Cubic operads}
If $\mathcal{A}ss=\Gamma(E)/(R_{\mathcal{A}ss})$ is the operad for associative algebras, the relations (\ref{cub})
 are the generating relations of $(R_{\mathcal{A}ss})(4)$. But these relations are following from the relations defining $(R_{\mathcal{A}ss})(3)=R_{\mathcal{A}ss}$. In Definition \ref{cubic associative}, we do not assume that the algebra is associative. It is clear that (\ref{cub})  do not implies associativity. From the relations (\ref{cub}) we can define a binary cubic operad $\mathcal{A}ss\mathcal{C}ubic$. 

Let $E$ be a $\K[\Sigma_2]$-module and $\Gamma(E)$ the free operad generated by $E$. Consider a $\K[\Sigma_4]$-submodule $R$ of $\Gamma(E)(4)$. Let $\mathcal{R}$ the ideal of $\Gamma(E)$ generated by $R$. We have
$$\mathcal{R}=\{\mathcal{R}(n), n \in \N^*\}$$
with $\mathcal{R}(1)=\{0\}, \ \mathcal{R}(2)=\{0\},  \ \mathcal{R}(3)=\{0\}, \ \mathcal{R}(4)=R.$
\begin{definition}
We call cubic operad generated by $E$ and defined by the relations $R\subset{\Gamma(E)(4)}$, the operad $\p(E,R)$ given by
$$\p(E,R)(n)=\displaystyle \frac{\Gamma(E)(n)}{\mathcal{R}(n)}.$$
\end{definition}
For example, the operad $\mathcal{A}ss\mathcal{C}ubic$ is the cubic operad generated by $E=\K[\Sigma_2]$ and the $\K[\Sigma_4]$-submodule of relations $R$  generated by the vectors
$$
\left\{
\begin{array}{l}
((x_1x_2)x_3)x_4-(x_1(x_2x_3))x_4,
(x_1(x_2x_3))x_4-x_1((x_2x_3)x_4),
x_1((x_2x_3)x_4)-x_1(x_2(x_3x_4)),\\
x_1(x_2(x_3x_4))-(x_1x_2)(x_3x_4),
(x_1x_2)(x_3x_4)-((x_1x_2)x_3)x_4.
\end{array}
\right.
$$
Thus we have $\mathcal{A}ss\mathcal{C}ubic(2)=\K[\Sigma_2]$, $\mathcal{A}ss\mathcal{C}ubic(3)=\K[\Sigma_3]$, and
$\mathcal{A}ss\mathcal{C}ubic(4)=\displaystyle \frac{\Gamma(E)(4)}{\mathcal{R}(4)}$
is the $24$-dimensional $\K$-vector space generated by $\{\left((x_{\sigma(1)}x_{\sigma(2)})x_{\sigma(3)}\right)x_{\sigma(4)}, \ \sigma \in  \Sigma_4\}$.
\medskip

\begin{definition}
The cubic operad $3\mathcal{N}ilp$ is defined by $3\mathcal{N}ilp(2)=sgn_2$, $3\mathcal{N}ilp(3)=\mathcal{L}ie(3)$ and $3\mathcal{N}ilp(4)=\{0\}$, where $\mathcal{L}ie$ is the operad for Lie algebras.
\end{definition}
In particular any $3\mathcal{N}ilp$-algebra is a $3$-step nilpotent Lie algebra. But this operad is cubic so we don't know direct links between an operadic cohomology and the cohomology of deformations.
\medskip

\subsection{Characteristic sequence $(3,\cdots,3,1)$}

Let $\g_{p,0,1}$ be the $(3p+1)$-dimensional nilpotent Lie algebra whose Lie bracket is
\begin{equation}
\label{clas331}
\left\{
\begin{array}{l}
[X_1,X_{3i-1}\rbrack=X_{3i}, \ [X_1,X_{3i}]=X_{3i+1}, \ i=1,\cdots,p .\\
\end{array}
\right.
\end{equation}
Then, any $(3p+1)$-dimensional $3$-step nilpotent Lie algebra $\g$ whose characteristic  sequence is $(3,\cdots,3,1),$ that is $\g \in \mathcal{F}^{3p+1,3}_{p,0,1},$ is isomorphic to a linear deformation of $\g_{p,0,1}$. If $\mu$ (resp. $\mu_0$)  is the Lie bracket of $\g$ (resp. of $\g_{p,0,1}$ ), then $\mu=\mu_0+t\varphi$ where the bilinear form $\varphi$ satisfies (\ref{3step}). Since  $\varphi$ has to satisfy $\delta_{C,\mu_0}^2(\varphi)=0$ and  $\delta_{R,\mu_0}^2(\varphi)=0$, this $2$-cocycle is cohomologous to 
the bilinear form, always denoted by $\varphi$:
\begin{equation}
\left\{
 \label{p01}
\begin{array}{l}
\ds
\medskip
\varphi(X_{3i-1},X_{3j})=\sum_{k=1}^pa_{3i-1,3j}^kX_{3k+1}, \ 1 \leq i \leq j \leq p,\\
\medskip
\ds \varphi(X_{3i},X_{3j-1})=\sum_{k=1}^p a_{3i,3j-1}X_{3k+1}, \ 1 \leq i < j \leq p,\\
\medskip
\ds \varphi(X_{3i-1},X_{3j-1})=\sum_{k=1}^p (a_{3i,3j-1}+a_{3i-1,3j})X_{3k}+\sum_{k=1}^p a_{3i-1,3j-1}X_{3k+1}, \\ 
\qquad \qquad \qquad \qquad  \qquad \qquad \qquad \qquad \qquad \qquad \qquad \ 1 \leq i < j \leq p.\\
\end{array}
\right.
\end{equation}
In particular, we deduce that 
$$\dim H^2_{CR}(\g_{p,0,1},\g_{p,0,1})=\ds \frac{p^2(3p-1)}{2}.$$
If $f \in GL(3p+1,\K)$, then the conditions $\delta_{C,\mu_0} f(X_1,X_i)=0$ imply that 
$$
\left\{
\begin{array}{l}
f(X_{3i})=[f(X_1),X_{3i-1}]+[X_1,f(X_{3i-1})], \\
f(X_{3i+1})=[f(X_1),X_{3i}]+[X_1,f(X_{3i})].
\end{array}
\right.
$$ Then 
$$\dim B^2_{CR}(\g,\g) \leq (3p+1)(p+1)$$
for any $\g \in \mathcal{F}^{3p+1,3}_{p,0,1}.$ Since any $2$-cocycle of $\g$ is also given by (\ref{p01}), we deduce
\begin{proposition}
If $p \geq 4$, any $\g \in \mathcal{F}^{3p+1,3}_{p,0,1}$ is not rigid. 
\end{proposition}
In particular, since $\mathcal{F}^{3p+1,3}_{p,0,1}$ is isomorphic to a linear space of dimension $ \frac{p^2(3p-1)}{2}$, this family is the orbit of a family with at least $\frac{3p^3-7p^2-8p-2}{2}$ parameters. 

Let us consider  the case $p=2,$ that is, of the dimension $7$. We have
$$
\left\{
\begin{array}{l}
\varphi(X_2,X_3)=a_1X_4+b_1X_7,\\
\varphi(X_3,X_5)=a_2X_4+b_2X_7,\\
\varphi(X_2,X_6)=a_3X_4+b_3X_7,\\
\varphi(X_5,X_6)=a_4X_4+b_4X_7,\\
\varphi(X_2,X_5)=(a_2+a_3)X_3+a_5X_4+(b_2+b_3)X_6+b_5X_7.\\
\end{array}
\right.
$$
In particular, computing $\delta_{C,\mu_0} f (X_i,X_j)$ with the condition $\varphi (X_1,Y)=0$ for any $Y$, we find
$$\dim H^2(\g_{2,0,1},\g_{2,0,1})=10$$
and the family $\mathcal{F}^{7,3}_{2,0,1}$ is reduced to a $10$-dimensional regular algebraic variety parametrized by the constants $(a_i,b_i), \ i=1,\cdots,5.$

\medskip

\begin{proposition}
The $7$-dimensional nilpotent Lie algebra $\g$ whose Lie bracket is
\begin{equation}
\label{rigid}
\left\{
\begin{array}{l}
[X_1,X_{3i-1}\rbrack=X_{3i}, \ [X_1,X_{3i}]=X_{3i+1}, \ i=1,2 ,\\
\lbrack X_2,X_3 \rbrack=X_4+X_7,\
\lbrack X_3,X_5 \rbrack=X_7,\\
\lbrack X_5,X_6 \rbrack=X_4,\\ 
\lbrack X_2,X_5\rbrack=X_6.\\
\end{array}
\right.
\end{equation}
is rigid in $3Nilp_7.$  Any $7$-dimensional Lie algebra of characteristic sequence $(3,3,1)$ is a contraction of $\g$.
\end{proposition}
In fact $\dim H^2(\g,\g) = 0 \footnote{During a talk given this summer, someone has pointed out to us that the example written in the article did not correspond to the one that we exposed. A term was missing. It is corrected in this version.}$  

We deduce that the classification of the Lie algebras of $\mathcal{F}^{7,3}_{2,0,1}$ is given by
$$
\begin{array}{l}
(a_1,a_2,a_3,a_4,a_5,b_1,b_2,b_3,b_4,b_5)  \in  \{  (0,0,0,0,0,0,0,0,0,0) ,
       (0,0,0,0,0,0,0,0,0,1), \\
       (0,0,1,0,0,0,0,0,0,0), (0,0,0,1,0,0,0,0,0,1),(0,0,0,1,0,0,0,0,0,0),\\
       (0,0,0,1,0,1,0,0,0,0),
       (1,0,0,0,0,1,0,0,1,0),(1,0,0,0,0,0,0,0,1,0),\\(1,0,0,0,1,0,1,0,0,0),
      (0,0,0,0,0,1,0,0,1,0), 
      (0,0,1,1,0,0,0,0,1,0),\\(1,0,0,1,0,0,1,0,0,0) , (1,0,0,1,0,0,0,0,0,1),(1,0,0,0,0,0,1,0,0,0),\\
      (0,0,0,0,0,0,0,0,1,0),(1,0,0,0,0,0,0,0,0,1) \}
 \end{array}
 $$

\subsection{Algebras attached to a $3$-step nilpotent Lie algebra}

If $\g$ is a $3$-step nilpotent Lie algebra, there exist $k_1,k_2,k_3$ such that $\g \in \mathcal{F}_{k_3,k_2,k_1}^{n,3} $ and it is a linear deformation of $\g_{k_3,k_2,k_1}$. Then the brackets $\mu$ of $\g$ and $\mu_0$ of  $\g_{k_3,k_2,k_1}$ satisfy $\mu=\mu_0+t\varphi$ and Relations (\ref{3step}). In particular we have
\begin{equation}
\label{alg}
(\mu_0 \circ_1\varphi) \circ_1 \varphi +( \varphi \circ_1 \varphi )\circ_1 \mu_0  +(\varphi\circ_1 \mu_0 )  \circ_1 \varphi =0.
\end{equation}
This relation can be interpreted as a definition of an algebra structure on $\K^n.$ This algebra is not necessarily a Lie algebra and it is associated with $\mu_0.$ 
\begin{definition}
Let $\g=(\K^n, \mu)$ be a $n$-dimensional $3$-step nilpotent Lie algebra. We call algebra attached to $\g$ any $n$-dimensional $\K$-algebra whose bilinear multiplication $\varphi$ satisfies the following identity:
$$(\mu \circ_1\varphi) \circ_1 \varphi + (\varphi \circ_1 \varphi )\circ_1 \mu  +(\varphi\circ_1 \mu )  \circ_1 \varphi =0.$$
\end {definition}
\noindent For example, if $\g=\g_{p,0,1}=(\K^{3p+1},\mu_0)$, the algebras whose multiplication is (\ref{p01}) is attached to $\g_{p,0,1}$.

For this type of algebras, we can define a cohomology of deformations which provides new invariants of $\g$. For this cohomology, $\mu_0$ is always a $2$-cocycle. 

\noindent{\bf Remark.} There is a class of cubic algebras which is  really interesting: the Jordan algebras. Recall that a $\K$-Jordan algebra is a commutative algebra satisfying the following identity
$$x(yx^2)=(xy)x^2.$$
Since $\K$ is of zero characteristic, linearizing this identity, we obtain
$$((x_2x_3)x_4)x_1+((x_3x_1)x_4)x_2+((x_1x_2)x_4)x_3-(x_2x_3)(x_4x_1)-(x_3x_1)(x_4x_2)-(x_1x_2)(x_4x_3)=0.$$
This relation is cubic. It is invariant by the permutations $Id, \tau_{12}, \tau_{13},\tau_{23}, c, c^2$ where $c$ is the cycle $(123).$  This permits to determinate an associated cubic operad called $\mathcal{J}ord$.  It is defined by
$$\mathcal{J}ord(2)=1\! \! 1, \ \mathcal{J}ord(3)=1\! \! 1 \oplus V, \ \mathcal{J}ord(4)=\displaystyle \frac{\Gamma(1\! \! 1)(4)}{\mathcal{R}(4)}$$
where $1\! \! 1$ is the identity representation of $\Sigma_2$ and  the $\K[\Sigma_4]$-module $\mathcal{R}(4)$ is generated by the vector
$$((x_2x_3)x_4)x_1+((x_3x_1)x_4)x_2+((x_1x_2)x_4)x_3-(x_2x_3)(x_4x_1)-(x_3x_1)(x_4x_2)-(x_1x_2)(x_4x_3)$$
is a vector space of dimension $4$.  Since $\dim \Gamma(1\! \! 1)(4)= 15$, thus
$\displaystyle \frac{\Gamma(1\! \! 1)(4)}{\mathcal{R}(4)}$ is of dimension $11$.
But, as in the $\mathcal{N}ilp_3$ operad case, we do not know a direct process giving an operadic cohomology or a link between this cohomology and the deformation cohomology. But we can directly define the space of bilinear forms associated with a linear deformation. This space corresponds to the $2$-cocycles in the deformation cohomology.  Using notations of \cite{Goze-Remm-Nonass}, if $\mu$ is the Jordan multiplication, then the Jordan identity is writen 
$$((\mu \circ_1 \mu )\circ_1\mu  - \mu \circ (\mu \otimes \mu) )\circ \Phi_v=0$$ where $v$ is the following vector of the group algebra $\K[\Sigma_4]$, $v= (2341)+ (3142) + \tau_{34}$. 
We deduce that a $2$-cocycle is defined by the identity
$$((\varphi \circ_1 \mu) \circ_1\mu+(\mu\circ_1 \varphi )\circ_1\mu+(\mu \circ_1 \mu )\circ_1\varphi- \mu \circ (\mu \otimes \varphi +\varphi \otimes \mu)+\varphi  \circ (\mu \otimes \mu)) \circ \Phi_v=0.$$

\noindent{\bf Remarks.}
\begin{enumerate}
\item We can generalize this process to define $(n-1)$-associative (binary) algebras: we consider the relations defining the $\Sigma_n$-module $\mathcal{A}ss(n)$ of the quadratic operad $\mathcal{A}ss$ and define, as above, an algebra with a multiplication which is a bilinear map $\mu$ (nonassociative), satisfying the previous relations where $\mu$ occurs $n-1$ times in each term of the relations. This algebra will be called $(n-1)$-associative (binary) algebra. If the Lie bracket of a algebra $\g$ is also $(n-1)$-associative, we prove in a similar way than for the cubic associative case that $\g$ is a nilpotent Lie algebra of nilindex $n-1$.

\item There exits another notion of associativity for $n$-ary algebras (an $n$-ary algebra is a vector space with  a multiplication which is an $n$-linear map), the total associativity. For example, a totally associative $3$-ary algebra has a ternary multiplication, denoted $xyz$, satisfying the relation:
    $$(xyz)tu=x(yzt)u=xy(ztu)$$ for any $x,y,z,t,u$.
  The corresponding operad is studied in \cite{Nico-Elis-JOA} and \cite{Elis-tartu}. Let $\g$ be a Lie algebra. We have the notion of Lie triple product given by $[[x,y],z].$ If we consider the vector space $\g$ provided with the $3$-ary product given by the Lie triple product, then $\g$ is a $3$-Lie algebra (\cite{Mic-Nic-Elis-Rabat}). Let us suppose now that the Lie triple bracket of $\g$ is a totally associative product. This implies
    $$
    [[[[X,Y],Z],T],U]=[[X,[[Y,Z],T]],U]=[[X,Y],[[Z,T],U]].$$
    But
    $$
    \begin{array}{lll}
    [[X,Y],[[Z,T],U]]&=&-[[[Z,T],U],[X,Y]]\\
    &=&[X,[Y,[[Z,T],U]]]+[Y,[[[Z,T],U],X]]\\
    &=& [[[[Z,T],U],Y],X]-[[[[Z,T],U],X],Y]\\
    &=&[[Z,[[T,U],Y]],X]-[[Z,[[T,U],X]],Y].
    \end{array}
    $$
    We deduce
    $$
    \begin{array}{lll}
    $$[[X,[[Y,Z],T]],U] &=&[[Z,[[T,U],Y]],X]-[[Z,[[T,U],X]],Y]\\
&=&2^5[[X,[[Y,Z],T]],U]-2^5[[X,[[Y,Z],U]],T]
    \end{array}
    $$
    Then
    $$(2^5-1)[[X,[[Y,Z],T]],U]=2^5[[X,[[Y,Z],U]],T].$$
    This implies
    $$[[X,[[Y,Z],T]],U]=0=[[[[X,Y],Z],T],U]=[[X,Y],[[Z,T],U]].$$
    The Lie algebra is $4$-step nilpotent.
\end{enumerate}


\begin{thebibliography}{99}

\bibitem{Ancochea-Goze} Ancoch\'ea-Berm\'udez, J.M; Goze, M. Sur la classification des alg\`ebres de Lie nilpotentes de dimension 7. C. R. Acad. Sci. Paris S\'er. I Math. 302 (1986), no. 17, 611–-613.

\bibitem{Ge} Gerstenhaber, Murray. On the deformation of rings and algebras. Ann. of Math. (2) 79 1964 59--103. 

\bibitem{Gong} Gong, Ming-Peng. Classification of nilpotent Lie algebras of dimension 7 (over algebraically closed fields and R). Thesis (Ph.D.). University of Waterloo (Canada). 1998. 165 pp. 

\bibitem{MGoran} Goze, M. Alg\`ebres de Lie r\'eelles et complexes, in {\it{Alg\`ebre, dynamique et analyse pour la g\'eom\'etrie: aspects r\'ecents}}. Editions Ellipses. 2010.

\bibitem{Mic-Nic-Elis-Rabat} Goze, M.; Goze, N.;  Remm, E. $n$-Lie algebras. African Journal of Maths and Physic 8. (2010). 17--28.

\bibitem{GK}  Goze, M; Khakimdjanov, Y. Nilpotent Lie algebras. Mathematics and its Applications, 361. Kluwer Academic Publishers Group, Dordrecht, 1996. xvi+336 pp. 

\bibitem{GP}  Goze, M; Piu, P. Une caract\'erisation riemannienne du groupe de Heisenberg.  Geom. Dedicata 50 (1994), no. 1, 27--36.

\bibitem{G.R1} Goze, M.;  Remm, E. Lie-admissible algebras and
operads. J. Algebra {\bf 273} (2004), no. 1, 129--152.

\bibitem{Goze-Remm-Nonass} Goze, M.;  Remm, E.  A class of nonassociative algebras.
Algebra colloquium, {\bf  14} no.2, (2007), 313--326.

\bibitem{Goze-Remm-Valued} Goze, M.;  Remm, E.  Valued deformations of algebras. J. Algebra Appl. 3 (2004), no. 4, 345--365.

\bibitem{GR-DGA} Goze, M.;  Remm, E.  Contact and Frobeniusian forms on Lie groups. Differential Geom. Appl. 35 (2014), 74--94.

\bibitem{Nico-Elis-JOA} Goze, N.; Remm, E. Dimension theorem for the free ternary partially associative algebras and applications. J. Algebra {\bf 348} (2011), no. 1, 14--36.


\bibitem{Markl} Markl, M. Deformation theory of algebras and their diagrams. CBMS Regional Conference Series in Mathematics, 116. Published for the Conference Board of the Mathematical Sciences, Washington, DC; by the American Mathematical Society, Providence, RI, 2012. x+129 pp.


\bibitem{M.S.S} Markl, M.; Shnider,  S.;  Stasheff, J. Operads in algebra,
topology and physics. Mathematical Surveys and Monographs, {\bf 96}. American
Mathematical Society, Providence, RI, 2002.

\bibitem{Elis-tartu} Remm, E. On the NonKoszulity of ternary partially associative Operads.
Proceedings of the Estonian Academy of Sciences, {\bf 59}, 4, (2010) 355--363.

\bibitem{NR} Nijenhuis, Albert; Richardson, R. W., Jr. Deformations of Lie algebra structures. J. Math. Mech. 17 1967 89--105.

\bibitem{Schaps}  Schaps, M. Deformations of algebras and Hochschild cohomology. {\it Perspective in ring theory, Antwerpen 1987}. 41--58. NATO adv. Sci. Ser. C. Math. Phys.Sci. 233. Kluwer Acad. Publish. DordrechT. 1988.


\bibitem{Zaili} Zaili, Y; Shaoqiang, D. The classification of two step nilpotent complex Lie algebras of dimension $8.$ Czechoslovak Math. J. 63(138) (2013), no. 3, 847--863.



 
\bibitem{Vergne} Vergne, Mich\`ele.
\newblock Cohomologie des alg\`ebres de Lie nilpotentes. Application  l'\'etude de la vari\'et\'e des alg\`ebres de Lie nilpotentes.
\newblock {\em C. R. Acad. Sci. Paris }, Sr. A-B 267 1968.






\end{thebibliography}
\end{document}